\documentclass[a4paper,12pt]{article}
\usepackage[utf8]{inputenc}
\usepackage[english]{babel}
\usepackage{amssymb}
\usepackage[nottoc]{tocbibind}
\usepackage{amsmath}
\usepackage{breqn}
\usepackage{authblk}
\usepackage[a4paper, total={6in, 8in}]{geometry}
\title{Synchronization Fronts in a Spatially Extended System of Hybrid Rayleigh-van der Pol Oscillators}

\author[1,*]{Carles Tardío Pi}
\author[2]{Jorge A. Castillo Medina}
\author[3]{Pablo Padilla Longoria}
\affil[1]{\footnotesize Laboratorio de Biología de Sistemas, Centro de Ciencias Genómicas, Universidad Nacional Autónoma de México, Av. Universidad s/n Col. Chamilpa 62210, Cuernavaca, Morelos, México.}
\affil[2]{\footnotesize Facultad de Matemáticas, Universidad Autónoma de Guerrero, Carlos E. Adame  54, Col.  Garita, C.P.39650, Acapulco, Guerrero, México.}
\affil[3]{\footnotesize Insituto de Investigaciones en Matemáticas Aplicadas y Sistemas, Universidad Nacional Autónoma de México, Cto. Escolar 3000, C.U., Coyoacán, 04510 Ciudad de México, México.}

\affil[*]{Corresponding author: carlestapi@gmail.com}
\date{\footnotesize \textit{Preprint} - 
May 6, 2022}

\usepackage{graphicx}
\usepackage{float}

\usepackage{caption}
\usepackage{subcaption}

\usepackage{listings,newtxtt}
\usepackage{listings}
\usepackage{xcolor}

\definecolor{codegreen}{rgb}{0,0.6,0}
\definecolor{codegray}{rgb}{0.5,0.5,0.5}
\definecolor{codepurple}{rgb}{0.58,0,0.82}
\definecolor{backcolour}{rgb}{0.95,0.95,0.92}
\lstdefinestyle{mystyle}{
    backgroundcolor=\color{backcolour},   
    commentstyle=\color{codegreen},
    keywordstyle=\color{magenta},
    numberstyle=\tiny\color{codegray},
    stringstyle=\color{codepurple},
    basicstyle=\ttfamily\footnotesize,
    breakatwhitespace=false,         
    breaklines=true,                 
    captionpos=b,                    
    keepspaces=true,                 
    numbers=left,                    
    numbersep=5pt,                  
    showspaces=false,                
    showstringspaces=false,
    showtabs=false,                  
    tabsize=2
}

\begin{document}

\maketitle

\section*{Abstract}

Numerous biological systems exhibit transitions to synchronised oscillations via a population-density-dependant mechanism known as quorum sensing. Here we propose a model system, based on spatially distributed limit-cycle
oscillators, that allows us to capture the dynamics of synchronization fronts by taking the continuum limit as the number of coupled oscillators tends to infinity. We explore analytically a family of wavefront-type solutions to the system in terms of model parameters and corroborate its validity via a numerical finite-difference method algorithm. Finally, we review our results in light of some experimental systems in synthetic biology based on a synchronised quorum of genetic clocks, coupled via a diffusive auto-inducer signalling molecule analogue to our Hooke-like coupling scheme within the proposed hybrid Rayleigh-van der Pol model.

\providecommand{\keywords}[1]
{
  \small	
  \textbf{\textit{Keywords---}} #1
}
\keywords{Coupled Oscillators, Genetic Clocks, Quorum Sensing, Synchronization Fronts, Limit-cycle, Travelling Waves, Autoinducer, Repressilator, Synthetic Biology.}

\section{Introduction}

Synchronization in a large population of interacting oscillatory elements is a remarkable phenomenon that has attracted attention in multiple research areas \cite{pikovsky2007synchronization}. One of the first historical reports regarding this phenomenon was made by a Dutch physician named Engelbert Kaempfer after a journey to Siam during the second half of the 17th century. Kaempfer describes an event which happens if a swarm of fireflies occupies a tree: the insects “hide their lights all at once, and a moment after make it appear again with the utmost regularity and exactness" \cite{bailey1988kaempfer}. This mechanism of synchronous emission of light pulses by a population of fireflies has been studied later on by biologists \cite{buck1968mechanism} and mathematically through the study of simplified models of coupled oscillators \cite{strogatz1993coupled}. 

A parallel biological phenomenon, but within the microscopic domain, has been investigated recently through the use of repressilator based synthetic gene regulatory networks of synchronized clocks coupled via quorum sensing \cite{danino2010synchronized}\cite{prindle2012sensing}\cite{mcmillen2002synchronizing}. These systems made of bacterial colonies are engineered to produce synchronised oscillations through a fluorescent protein reporter within a microfluidic device.

In both phenomena, composed by a group of interacting biological clocks and their environment, frequency synchronization is a logical approach since we are interested in states for which the temporal rates of change of the phases are equal, rather than the phases themselves.

One of the main models to study those kinds of systems that exhibit synchronization was first introduced by Winfree \cite{winfree1967biological}, and later refined by Kuramoto and others \cite{1985kuramoto}\cite{hong2002synchronization}\cite{hong2005collective}. (A comprehensive review and perspectives of the field can be found in \cite{pikovsky2015dynamics}). In the existing literature systems of fully coupled oscillators have mostly been considered due to its analytical tractability and simplicity . On the other hand, systems of spatially extended oscillators with local couplings have received less attention even though they are more realistic and frequently observed in nature \cite{hong2004collective}\cite{daido1988lower}\cite{bahiana1994order}\cite{aoyagi1991frequency}\cite{mirollo1990amplitude}.

Here we propose a system based on spatially distributed limit-cycle oscillators, which have been extensively used as a model of physical and biological clocks, based on a hybrid Rayleigh-van der Pol (RvdP) oscillator. The main idea is to capture the dynamics of synchronization fronts by taking the continuum limit as the number of coupled oscillators tends to infinity. This approach
makes sense in many biological systems such as colonies of bacteria in which a network-like coupling is generated via a quorum sensing mechanism. Usually when the population density exceeds a certain threshold then the corresponding local levels of signal can induce an spatially extended synchronised response in gene expression throughout the entire population \cite{li2012quorum}.

\section{Travelling Wave Solution to a Spatially Extended Rayleigh-van der Pol Equation}

We firstly consider the second order homogeneous van der Pol differential equation, describing an oscillator with non-linear damping,

\begin{equation}\label{vdppol0}
 \ddot{u}+\omega_0^2u-\varepsilon(1-u^2)\dot{u}=0,
\end{equation}

\noindent where $u\in\mathbb{R}$ is the dynamical variable, $\omega_0$ is the linear frequency and $\varepsilon>0$ is a non-linear scalar parameter of the system. An important aspect of this equation is the relation between the negative damping proportional to $-\dot{u}$, and the non-linear damping, $u^2\dot{u}$, which leads to the existence of self-sustained limit cycle oscillations.

A simple way to understand this equation is given by applying a transformation into a rotating reference frame \cite{lee2013quantum}, via 

\begin{equation}\label{vdp1}
    u(t)=\sqrt{2}\left( \alpha(t)e^{i\omega_0t}+\alpha^*(t)e^{-i\omega_0t}\right),
\end{equation}

\noindent with a complex amplitude that changes slowly as

\begin{equation}
    \alpha=\sqrt{2}(u+iv),
\end{equation}

\noindent where  $u$ and $v=\dot{u}$ denote the position and the conjugate moment respectively. Within this rotating frame one can neglect the terms corresponding to fast oscillations in (\ref{vdp1}) if $\varepsilon <<1$ holds. This allows us to write equation (\ref{vdp1}) as a generalised Stuart-Landau type equation 

\begin{equation}
    \dot{\alpha}(t)=\frac{\varepsilon}{2}\left(1-|\alpha(t)|^2 \right)\alpha(t),
\end{equation}

\noindent which describes the oscillator dynamics. Looking for  stationary states, i.e. when  $\dot{\alpha}(t)=0$, this equation admits a limit cycle defined by $|\alpha(t)|^2=1$.

In our case we introduce a modified version of equation (\ref{vdppol0}), namely

\begin{equation}\label{equation2}
    \ddot{u}+\omega_0 ^2 u-\varepsilon(1- u^2)\dot{u} - \delta\dot{u}^3 =0,
\end{equation}

\noindent which is of the form $\ddot{u}=f(u,\dot{u})$.

First of all we can observe that equation (\ref{equation2}) is a combination of the previous van der Pol equation  (\ref{vdp1}) and the classic Rayleigh equation defined by

\begin{equation}\label{rayleigh}
    \ddot{u}+\omega_0^2u- \varepsilon(1-\dot{u}^2)\dot{u}=0.
\end{equation}

In that sense the proposed equation (\ref{equation2}) contains the cubic term proportional to $\dot{u}^3$ which is contained in the Rayleigh equation (\ref{rayleigh}) and not in the van der Pol equation (\ref{equation2}). Concurrently it contains the  $u^2\dot{u}$ term which the Rayleigh equation (\ref{rayleigh}) doesn't have, but the van der Pol equation (\ref{equation2}) does. This hybrid Rayleigh-van der Pol equation has been studied in  an approximate way within the context of bipedal motion modelling \cite{haken1985theoretical}.

\begin{figure}[H]
    \centering
    \includegraphics[width=0.75\textwidth]{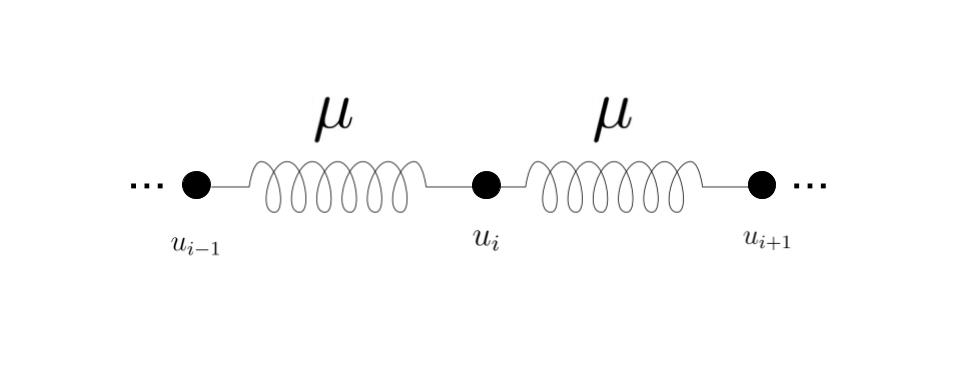}
    \caption{System diagram where each oscillator is coupled to its first neighbours via an elastic coupling constant.
    }
    \label{muellesvdp}
\end{figure}

Here we analyse equation (\ref{equation2}) starting from the case in which $\delta=1/\omega_0^2$ and considering the spatial location of each oscillator explicitly. To begin with, we consider a discrete system assuming that each element $u_i$ is coupled to its nearest neighbours  $u_{i-1}$ and $u_{i+1}$ via a coupling constant $\mu$ (see figure \ref{muellesvdp}). According to this scheme we obtain 
\begin{equation}\label{1scheme}
    \ddot{u}_i =f(u_i,\dot{u}_i)+\mu (u_{i+1}+u_{i-1}).
\end{equation}
Adding and subtracting the term $2u_i$
\begin{equation}
    \ddot{u}_i =f(u_i,\dot{u}_i)+\mu (u_{i+1}+u_{i-1}+2u_i-2u_i),
\end{equation}
we observe that the finite difference approximation of the second spatial derivative naturally appears.

In other words the term $u_{i+1}+u_{i-1}-2u_i$ is the discretisation of the partial derivative 
$$\frac{\partial ^2 u}{\partial x^2}\equiv u_{xx},$$
so the ODE becomes a PDE in the folling way,

\begin{equation}\label{odepde}
    u_{tt}\cong f(u,u_t)+2\mu u+ \mu u_{xx},
\end{equation}
and finally for one spatial dimension we obtain

\begin{equation}
    u_{tt}=-\omega_0^2 u + \varepsilon(1- u^2 )u_t- \delta u_t^3 + 2\mu u +\mu u_{xx},
\end{equation}
which can be rewritten as

\begin{equation}\label{eq1}
    u_{tt}=-\omega_0^2 u + \varepsilon(1- u^2- \delta u_t^2  )u_t+ 2\mu u +\mu u_{xx}.
\end{equation}

Now we show that a travelling wave of the form  \\
$u=sin(\xi)$, $\xi \equiv \omega_0 t+kx $ can be a solution of (\ref{eq1}) for an appropriate choice of the parameter $k$.

So we compute the time derivatives

\begin{equation}
    u_t=\omega_0 cos(\xi), u_{tt}=-\omega_0 ^2 sin(\xi),
\end{equation}

and the spatial ones as,

\begin{equation}
    u_x=a\, cos(\xi),u_{xx}=-k^2 sin(\xi).
\end{equation}

Then if we substitute in (\ref{eq1}) we get 

\begin{equation}\label{substitute}
    -\omega_0^2 sin\xi= -\omega_0^2 sin\xi+\varepsilon\left[1- \left(sin^2\xi + \frac{\omega_0^2cos^2\xi}{\omega_0^2}\right)\right]\omega_0 cos\xi+2\mu sin\xi+\mu(-k^2sin\xi),
\end{equation}
to which we obtain two solutions: the trivial one if $\mu=0$ and another one with wave number $k=\sqrt{2}$ ,

\begin{equation}\label{solucion}
u=sin(\omega_0 t+\sqrt{2}x),
\end{equation}

which is a solution to our system independently of the value of  $\varepsilon$.

\begin{figure}[H]
    \centering
    \includegraphics[width=0.75\textwidth]{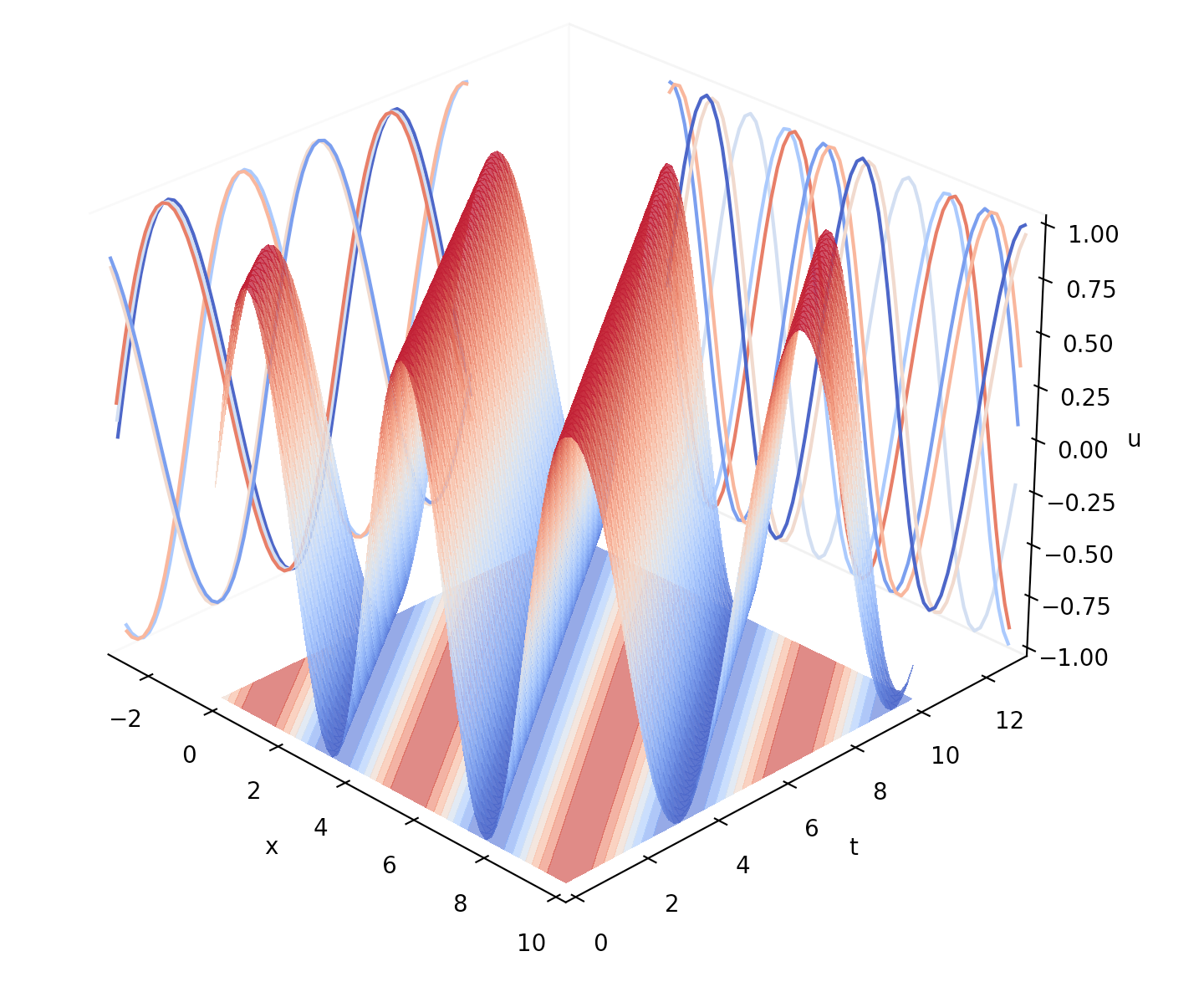}
    \caption{Analytic solution and orthogonal projections representing a travelling wave as in equation (\ref{solucion}), with parameters $k=\sqrt{2},\,\,\, \omega_0=1$.}
    \label{solucionanalitica}
\end{figure}

Next we can analyse the case where we take into consideration a 3-dimensional version of the previous coupling scheme of the system. By using a lattice coupling instead of the 1-dimensional represented in figure (\ref{muellesvdp}), the equation (\ref{1scheme}) stands as

\begin{equation}
    u_{i,j,k}=f(u_{i,j,k},\dot{u}_{i,j,k}))+\mu(u_{i+1}+u_{i-1})+\mu(u_{j+1}+u_{j-1})+\mu(u_{k+1}+u_{k-1}),
\end{equation}

as done previously in the 1-dimensional case we can add and subtract the $2\mu u_{i,j,k}$ term which gives us

\begin{equation}
u_{i,j,k}= f(u_{i,j,k},\dot{u}_{i,j,k}))+\mu(u_{i+1}+u_{i-1}+2u_{i}-2u_{i})+\mu(u_{j+1}+u_{j-1}+2u_{j}-2u_{j}) \\
+\mu(u_{k+1}+u_{k-1}+2u_{k}-2u_{k}),
\end{equation}

so as seen previously in equation (\ref{odepde}) the ODE becomes a PDE as 

\begin{equation}\label{laplace}
  u_{tt}\cong f(u,u_t)+6\mu u+ \mu( u_{xx}+ u_{yy}+ u_{yy})\cong f(u,u_t)+6\mu u+ \mu\nabla^{2}u,   
\end{equation}
being $\nabla^{2}$ the Laplacian operator defined as $\nabla^{2}=\frac{\partial^2 }{\partial x^2}+\frac{\partial^2 }{\partial y^2}+\frac{\partial^2 }{\partial z^2}$ .

Similarly to equation (\ref{substitute}), we find that a solution in the form $u=sin(\xi), \xi=\omega_0 t+k x $ in equation (\ref{laplace}) leads to $6\mu sin(\xi)+\mu(-3k^2sin(\xi))=0$, so we recover the exact same solutions $\mu=0$ and $k=\sqrt{6/3}=\sqrt{2}$ as in equation (\ref{solucion}).

Moreover, we can also analyse the spatially homogeneous case in equation (\ref{eq1}) by dropping the second spatial derivative in its last term, which leads to

\begin{equation}
    u_{tt}=-(\omega_0^2 -2\mu)u + \varepsilon(1- u^2- \delta u_t^2  )u_t.
\end{equation}

If we propose again an analytic solution in the form of a travelling wave but without the spatial term by $u=sin(\Omega t)$ we get

\begin{equation}
    -\Omega^2 sin(\Omega t)=-(\omega_0-2\mu)sin(\Omega t)+\varepsilon(1- sin^2(\Omega t)- \delta \Omega^2 cos(\Omega t)  ),
\end{equation}

which corresponds to the case where $\delta=1/\Omega^2$ and $\Omega=\sqrt{\omega_0-2\mu}$ independently of the value of the parameter $\varepsilon$.

Finally, with the aim of finding an asymptotic frequency we linearize equation (\ref{equation2}) and we take $\varepsilon<<1$ and $\delta<<1$ so we can neglect the dissipative term. Similarly, we can add to equation (\ref{equation2}) a coupling between its elements via a Laplacian matrix $L=[L_{ij}:i,j=1,...,n]$ for the graph  $\mathcal{G}$. In that sense, $L_{ij}=a_{ij}$ if $i\neq j$ and  $L_{ii}=\sum_{j=1}^{N} a_{ij}$. In this way we obtain the following equation 

\begin{equation}
    \ddot{u}=-\omega ^2 u+\mu \sum_{j=1}^{N}a_{ij}x_j, \quad i=1,...,N,
\end{equation}

which can be rewritten as

\begin{equation}
    \ddot{u}=-\left(\omega ^2 -\mu A_i\right)u_i,
\end{equation}

where $A_i=\sum_{j=1}^{N} a_{ij}$ for $i=1,...,N$. In the case that we have $a_{ij}=1$, corresponding to a completely connected network of oscillators, the system can be written as 

\begin{equation}
    \ddot{u}=-\left(\omega ^2 -\mu (N-1)\right)u.
\end{equation}

Supposing a synchronised regime has been established, we can have a synchronization frequency determined by

\begin{equation}
    \Omega_{sync}=-\sqrt{\omega ^2 -\mu (N-1)} .
\end{equation}

Apart from this, this reduction suggests that the synchronised oscillatory behaviour happens for sufficiently small values of $\mu$ and also when

\begin{equation}
    \omega ^2 -\mu (N-1)>0.
\end{equation}

If this wasn't the case we would have an exponentially growing behaviour in our system. Otherwise, we can observe that unless all $a_{ij}$ are equal we cannot extract these conclusions. So, to analyse the general case where not all $a_{ij}$ are equal, we could use for instance certain methods that introduce perturbations to the Fokker-Planck equation \cite{franci2018synchronization}. Nonetheless, such analysis is outside the scope of this paper.

In the next section we develop a numerical method based on a finite differences algorithm techniques in order to numerically check the analytic solution that we have found. Finally we interpret those results in the context of the synchronised quorum of genetic clocks experiments in synthetic biology \cite{danino2010synchronized}.

\section{Numerical Methods based on a Finite Difference Algorithm}

In order to find a solution to our system we develop a numerical simulation based on a finite differences method \cite{langtangen1999computational}. According tho this we start with the generic PDE 

\begin{equation}
    u_{tt}=\nabla\cdot(c^2\nabla u)+f,
\end{equation}

which is then transformed to lead to our equation (\ref{eq1}). We begin by discretizing our temporal domains 
$[0,T]$ and the spatial ones $[0,L]$ so

\begin{equation}
0=t_{0}<t_{1}<t_{2}<\cdots<t_{N_{t}-1}<t_{N_{t}}=T
\end{equation}
and 
\begin{equation}
0=x_{0}<x_{1}<x_{2}<\cdots<x_{N_{x}-1}<x_{N_{x}}=L.
\end{equation}

In this manner we can introduce for a distributed set of points the spacing constants $\Delta t$ y $\Delta x$ within a mesh to get 

\begin{equation}
x_{i}=i \Delta x, \quad i=0, \ldots, N_{x}, \quad t_{n}=n \Delta t,\quad n=0, \ldots, N_{t}.
\end{equation}

Then the solution $u(x,t)$ is obtained for the points in the mesh, were we introduce the notation $u_i^n$,
which approximates the exact solution in the points $\left(x_{i}, t_{n}\right)$ for $i=0, \ldots, N_{x}$ y $n=0, \ldots, N_{t}$.

The next step is to change the derivatives by finite differences according to 

\begin{equation}
\frac{\partial^{2}}{\partial t^{2}} u\left(x_{i}, t_{n}\right) \approx \frac{u_{i}^{n+1}-2 u_{i}^{n}+u_{i}^{n-1}}{\Delta t^{2}},
\end{equation}

\begin{equation}
\frac{\partial^{}}{\partial t^{}} u\left(x_{i}, t_{n}\right) \approx \frac{u_{i}^{n+1}-u_{i}^{n}}{\Delta t^{}},
\end{equation}

and we introduce the operator notation for the finite differences as

\begin{equation}
\left[D_{t} D_{t} u\right]_{i}^{n}=\frac{u_{i}^{n+1}-2 u_{i}^{n}+u_{i}^{n-1}}{\Delta t^{2}},
\end{equation}

\begin{equation}
\left[D_{t} u\right]_{i}^{n}=\frac{u_{i}^{n+1}- u_{i}^{n}}{\Delta t^{2}},
\end{equation}

and for the spatial case

\begin{equation}
\frac{\partial^{2}}{\partial x^{2}} u\left(x_{i}, t_{n}\right) \approx \frac{u_{i+1}^{n}-2 u_{i}^{n}+u_{i-1}^{n}}{\Delta x^{2}}=\left[D_{x} D_{x} u\right]_{i}^{n}.
\end{equation}

We do the same with first order time derivatives obtaining the algebraic version of our PDE from equation (\ref{eq1}) as
 
\begin{equation}
\left[D_{t} D_{t} u=(-\omega_0^{2}+2\mu)u+\left[
1-\left(u^2+\frac{1}{\omega_0^{2}}[D_t u]^2\right)\right][D_t u]+\mu[D_{x} D_{x}u]\right]_i^n.
\end{equation}

\begin{figure}[H]
    \centering
    \includegraphics[width=0.5\textwidth]{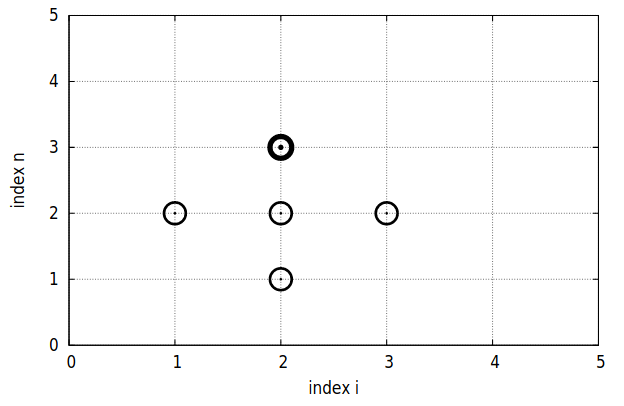}
    \caption{Spatio-temporal mesh. The circles show the connected points in a finite difference equation.
    }
    \label{egeo}
\end{figure}

That equation can be interpreted as a moving stencil which involves the values in the four pointsfrom its neighbourhood $u_{i}^{n+1}, u_{i \pm 1}^{n}, u_{i}^{n},$ y $u_{i}^{n-1}$, as we can see in the figure \ref{egeo}. 

In the same manner we have discretised the derivatives we need to discretize the initial conditions,

\begin{equation}
\frac{\partial}{\partial t} u\left(x_{i}, t_{n}\right) \approx \frac{u_{i}^{1}-u_{i}^{-1}}{2 \Delta t}=\left[D_{2 t} u\right]_{i}^{0},
\end{equation}

that can be written in algebraic notation as

\begin{equation}
\left[D_{2 t} u\right]_{i}^{n}=0, \quad n=0,
\end{equation}

leading to

\begin{equation} \label{init}
u_{i}^{n-1}=u_{i}^{n+1}, \quad i=0, \dots, N_{x}, n=0 ,
\end{equation}

and for the other initial condition

\begin{equation}
u_{i}^{0}=I\left(x_{i}\right), \quad i=0, \dots, N_{x}.
\end{equation}

Finally we propose a recursive algorithm where $u_i^{n+1}$ equals to the sum of the following $A,...,J$ terms as:

\begin{multline}\label{grande}
\underbrace{\left(\Delta t^2(-\omega_0^2+2\mu)+\Delta t-\frac{\Delta t^2\mu}{\Delta x^2}\right) }_{A}  u_i^n 
\underbrace{-\Delta t}_{B}(u_i^n)^2 
\underbrace{-\frac{1}{\omega_0^2\Delta t}}_{C}(u_i^n)^3 
\underbrace{-\Delta t}_{D}(u_i^n) 
\underbrace{-\Delta t}_{E}(u_i^n u_i^{n-1}) \\
\underbrace{-\frac{1}{\omega_0^2\Delta t}}_{F}(u_i^{n-1})^3 
\underbrace{-\frac{3}{\omega_0^2\Delta t}}_{G}(u_i^{n-1})(u_i^n)^3 
\underbrace{-\frac{3}{\omega_0^2\Delta t}}_{H}(u_i^{n+1})^2(u_i^n) 
\underbrace{+\frac{\mu\Delta t^2}{\Delta x^2}}_{I}(u_{i+1}^n) 
\underbrace{-\frac{\mu\Delta t^2}{\Delta x^2}}_{J}(u_{i-1}^n) ,
\end{multline}

where $A,...,J$ parameters are determined in our code (see section \ref{numerical}).

For the equation (\ref{grande}) we have a problem when $n=0$ since the formula for $u_{i}^{1}$ involves $u_{i}^{-1}$, which is not defined within the mesh. However, we can use the initial condition in (\ref{init}) for $n=0$ to eliminate $u_{i}^{-1}$.

Finally our algorithm can be summarised as:

\begin{enumerate}
    \item Calculate $u_{i}^{0}=I\left(x_{i}\right)$ for $i=0, \ldots, N_{x}$.
    \item Calculate $u_i^1$ from (\ref{grande}) and take $u_i^1=0$ in the boundary points $i=0$ e $i=N_{x}$ for $n=1,2, \dots, N-1$
    \item For each time step $n=1,2, \dots, N_{t}-1$ apply (\ref{grande}) to find $u_{i}^{n+1}$ for $i=1, \ldots, N_{x}-1$ and take $u_{i}^{n+1}=0$ for the boundary points $i=0, i=N_{x}$.
\end{enumerate}

\begin{figure}[H]
     \centering
     \begin{subfigure}[b]{\textwidth}
         \centering
         \includegraphics[width=\textwidth]{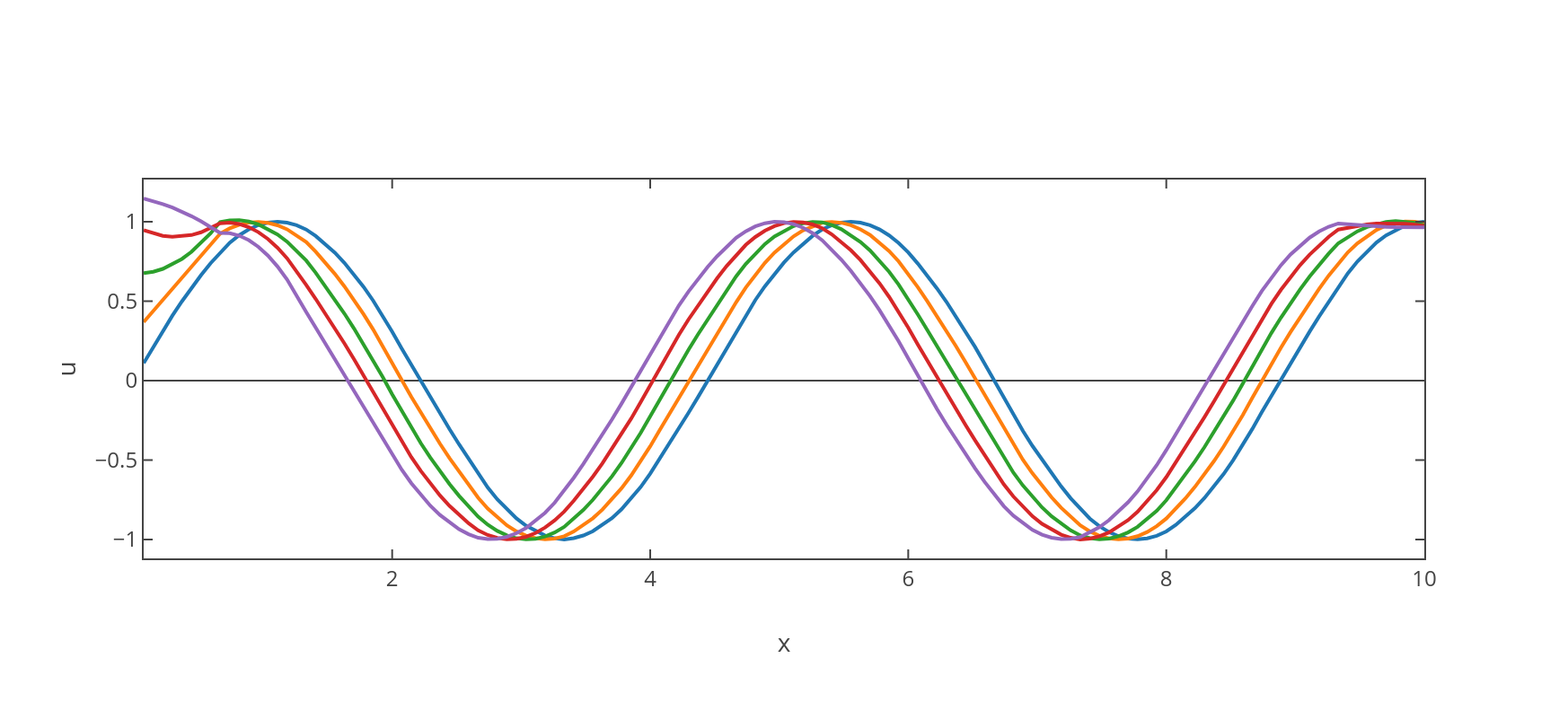}
     \end{subfigure}
     \hfill
     \begin{subfigure}[b]{\textwidth}
         \centering
         \includegraphics[width=\textwidth]{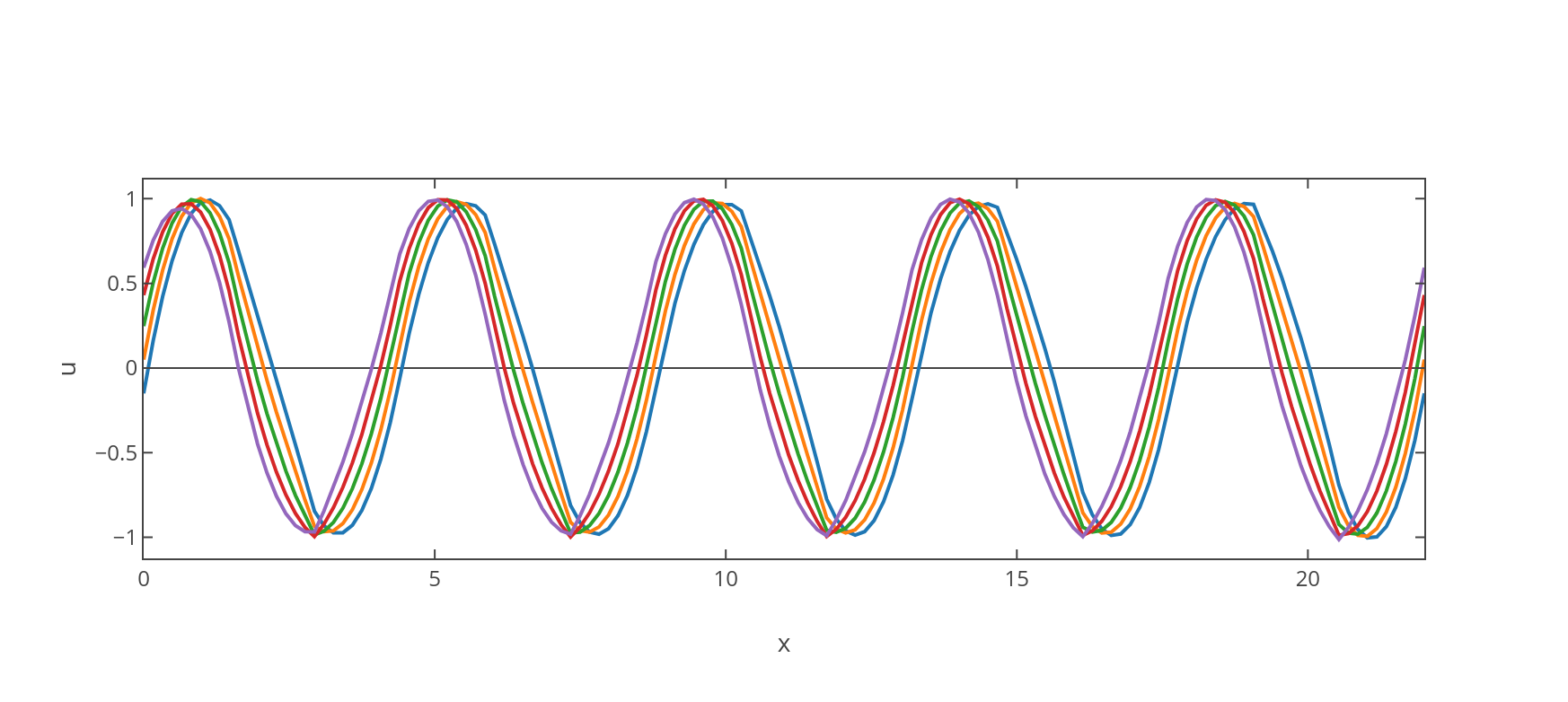}
     \end{subfigure}
    
     \caption{Numerical solution to equation (\ref{eq1}) plotted along 5 temporal steps every $\Delta t=0.2s$ (from light blue to purple), with parameters $\mu=\sqrt{2}$, $\varepsilon=1$, $\omega_0=2$ and $\delta=1/4$. On top figure, we select a $L=10$ spatial domain with no boundary conditions, and in the bottom we choose a $L=22$ domain with periodic boundary conditions such $u(x=0,t)=u(x=L,t)$.
  }
     \label{plotwave}
\end{figure}

This algorithm consists essentially in a moving stencil of finite difference elements through all the points in the mesh. Next we show the plotted results in Figure (\ref{plotwave}).

\section{Interpretation of the Travelling Wave Solution}

The synchronised quorum of genetic clocks synthetic system, which exhibits a range of oscillatory behaviours after a critical density of bacteria is reached, can be mainly framed by a network architecture that consists of an activator that activates its own repressor. Such motifs have been found to exist as regulatory modules in many genetic oscillators and circadian clock networks, studied both \textit{in vivo} and by mathematical models \cite{garcia2004modeling}.

Here we will focus on the compared analysis between our model, based on a system of spatially extended coupled oscillators, and some experimental realisations of the synchronised quorum system, such as  \cite{danino2010synchronized}\cite{prindle2012sensing}.  These kinds of  experimental designs make use of the previous repressilator-like architecture, additionally coupled to a system that synthesises a molecule known as an autoinducer (AI). This AI has the property that can diffuse through the cell membrane and allows cell-to-cell communication across the bacterial colony, in analogy to the Hooke-like coupling mechanism mediated by the parameter $\mu$ in our model.

The general solution that we found in (\ref{solucion}) corresponds to a travelling wave, which can be described as a periodic function in one-dimensional space that moves with constant phase-velocity determined by the dispersion relation $v_{\mathrm{p}}=\omega_0 /\sqrt{2}$, where $\omega_0$ corresponds to the angular frequency parameter selected within equation (\ref{eq1}). This fact has the implication that the wave travels undistorted since $v_{\mathrm{p}}=v_{\mathrm{g}}$, being the group velocity defined by $v_{\mathrm{g}} \equiv \partial \omega /\partial k$.

The solution which was obtained by imposing a coupling mechanism in between elements within the discretised version of the hybrid Rayleigh-van der Pol equation (\ref{equation2}) and then taking the continuous limit, does not depend on the coupling constant $\mu$. This might seem a counter-intuitive property of the system, nonetheless it seems to be a consequence of the robustness of the limit cycle under perturbations. Moreover, it makes sense from an evolutionary perspective, where one expects that the speed of propagation of a certain signal does not depend on the detailed features of the coupling, thus providing a structural advantage.
Furthermore, we can also notice that the solution in equation (\ref{eq1}) is also independent of the parameter $\varepsilon$, meaning that the limit cycle is unconstrained by the damping mechanism.

In addition, we can observe that switching the unitary term in equation (\ref{eq1}) by an extra parameter $A$, such that  

\begin{equation}\label{eq3}
    u_{tt}=-\omega_0^2 u + \varepsilon\left(A- u^2- \frac{u_t^2}{\omega_0^2}\right)u_t+ 2\mu u +\mu u_{xx},
\end{equation}

also allows for a more general travelling wave solution in the form of $u=Asin(\xi)$, $\xi \equiv \omega_0 t+\sqrt{2}x $, besides the trivial solutions $\mu=0$ and $A=0$. As a consequence, by selecting the parameters $\omega_0$ and $A$ within equation (\ref{eq3}), one can recover a family of harmonic waves with varying angular frequencies and maximum amplitudes constrained to a fixed wavenumber. 

Similarly, a wide range of spatio-temporal dynamics of the synchronised oscillators can be accomplished by varying some experimental parameters, such as the diffusion of the AI molecule, which turns out to modulate the amplitude and the period of the wave. Thus, for a high dissipation rate of AI both the amplitude and period increase. However, this result which is consistent with ‘degrade-and-fire’ oscillations observed previously in  intracellular oscillators \cite{stricker2008fast}, is not recovered in our model since its amplitude is ultimately decoupled from the frequency.

Additionally, we can note that the experimental results show a transition from an unsynchronised to a synchronised regime as the strength of coupling increases, caused by an increase in cell density. This characteristic property is recovered  in our model since the conditions for the emergence of the limit cycle are fulfilled once we take the continuum limit, allowing the number of oscillators to become infinite while the coupling remains local.

Finally, since the speed of wave front propagation in our system does not depend on the coupling parameter $\mu$, in contrast to the experimental case where $v\sim \sqrt{D}$ (being $D$ related to the diffusive coefficient of AI), a subsequent work will focus on studying the synchronization of other diffusive-like coupling mechanisms, where synchronization is a consequence of the dissipation in the coupling along with cases where there is an interaction between the inherent damping in the system and the coupling. 

\section{Conclusions}

In relation to the topics discussed in this paper, there is an issue regarding the emergence of the harmonic travelling wavefront solution that may be further clarified. The shape of the dispersion relation that we found, which turns out to be independent of the coupling proportionality constant between oscillators and also of the damping coefficient, may allude to an underlying robustness mechanism that needs to be interpreted. A possible way to do so would be in the context of biological soliton-like structures, defined as self-reinforcing solitary waves that travel at constant speed and without changing its shape \cite{kuwayama2013biological}. One procedure to asses the existence of this kind of structures would be to check out if our solution does not obey the superposition principle, which implies that the travelling wave propagates unperturbed in the presence of other wave structures.

In general terms, whilst analytically solvable examples of coupled-oscillator models may be limited in each biological topic under study, the continuum-limit picture itself, i.e the one-oscillator extended to a collective ensemble preserving the local coupling in a self-consistent manner, may remain useful for qualitative understanding the spatio-temporal dynamics of signal or energy propagation within large groups of coupled oscillators.

\section*{Supplementary Info.: Numerical Methods Code}\label{numerical}

\begin{lstlisting}[language=Python]
#!/usr/bin/env python3
# -*- coding: utf-8 -*-
"""
Title: FDM solution to RvdP equation
(Rev. with vectorized & Dirichlet, Neumann boundary conditions)
"""
import sys
import math
import numpy as np
import matplotlib.pyplot as plt

###

def solver(omega, dt, dx, mu, L, T, U_0, U_L, user_action=None,
           version='vectorized2'):
    """"""
    Nt = int(round(T/dt))
    t = np.linspace(0, Nt*dt, Nt+1)   # Mesh points in time

    Nx = int(round(L/dx))
    x = np.linspace(0, L, Nx+1)       # Mesh points in space

    # Make sure dx and dt are compatible with x and t
    dx = x[1] - x[0]
    dt = t[1] - t[0]

    if U_0 is not None:
        if isinstance(U_0, (float,int)) and U_0 == 0:
            U_0 = lambda t: 0
        # else: U_0(t) is a function
    if U_L is not None:
        if isinstance(U_L, (float,int)) and U_L == 0:
            U_L = lambda t: 0
        # else: U_L(t) is a function

    A = -(omega**2)*(dt**2)+dt+2+2*(dt**2)*mu-2*mu*(dt**2)/(dx**2)
    B = 2-dt
    C = dt*(1/(omega**2)-1)
    D = -dt*(1-3/(omega**2))
    E = 3*dt/(omega**2)
    F = -dt/(omega**2)
    G = mu*(dt**2)/(dx**2)
    H = mu*(dt**2)/(dx**2)

    u     = np.zeros(Nx+1)   # Solution array at new time level
    u_n   = np.zeros(Nx+1)   # Solution at 1 time level back
    u_nm1 = np.zeros(Nx+1)   # Solution at 2 time levels back

    Ix = range(0, Nx+1)
    It = range(0, Nt+1)

    import time;
    t0 = time.perf_counter()  # Measure CPU time

    # Solucion Analitica
    def u_exact(x, t):
        return math.sin(omega*t+math.sqrt(2)*x)

    # Load initial condition into u_n
    def I(x): return math.sin(math.sqrt(2)*x)
    #I = lambda x: u_exact(x, 0)

    # Load velocity initial condition 
    def V(x): return -0.5*omega*math.cos(math.sqrt(2)*x)

    for i in Ix:
        u_n[i] = I(x[i])

    if user_action is not None:
        user_action(u_n, x, t, 0)

    # Special formula for first time step

    for i in Ix[1:-1]:
        u[i] = A*u_n[i]+B*(u_n[i]-2*dt*V(x[i]))+C*(u_n[i]**3)+D*(u_n[i]**2)*(u_n[i]-2*dt*V(x[i]))\
        +E*(u_n[i])*((u_n[i]-2*dt*V(x[i]))**2)+\
        F*((u_n[i]-2*dt*V(x[i]))**3)+G*u_n[i+1]+H*u_n[i-1]

    i = Ix[0]
    if U_0 is None:
        # Set boundary values du/dn = 0
        # x=0: i-1 -> i+1 since u[i-1]=u[i+1]
        # x=L: i+1 -> i-1 since u[i+1]=u[i-1])
        ip1 = i+1
        im1 = ip1  # i-1 -> i+1
        u[i] = A*u_n[i]+B*(u_n[i]-2*dt*V(x[i]))+C*(u_n[i]**3)+D*(u_n[i]**2)*(u_n[i]-2*dt*V(x[i]))\
        +E*(u_n[i])*((u_n[i]-2*dt*V(x[i]))**2)+\
        F*((u_n[i]-2*dt*V(x[i]))**3)+G*u_n[ip1]+H*u_n[im1]
    else:
        u[0] = U_0(dt)

    i = Ix[-1]
    if U_L is None:
        im1 = i-1
        ip1 = im1  # i+1 -> i-1
        u[i] = A*u_n[i]+B*(u_n[i]-2*dt*V(x[i]))+C*(u_n[i]**3)+D*(u_n[i]**2)*(u_n[i]-2*dt*V(x[i]))\
        +E*(u_n[i])*((u_n[i]-2*dt*V(x[i]))**2)+\
        F*((u_n[i]-2*dt*V(x[i]))**3)+G*u_n[ip1]+H*u_n[im1]
    else:
        u[i] = U_L(dt)

    if user_action is not None:
        user_action(u, x, t, 1)

    # Update data structures for next step
    #u_nm1[:] = u_n;  u_n[:] = u  # safe, but slower
    u_nm1, u_n, u = u_n, u, u_nm1

    for n in It[1:-1]:
        # Update all inner points at time t[n+1]
        if version == 'scalar':
            for i in Ix[1:-1]:
                u[i] = A*u_n[i]+B*u_nm1[i]+C*((u_n[i])**3)+D*((u_n[i])**2)*(u_nm1[i])+E*(u_n[i])*((u_nm1[i])**2)+\
                F*((u_nm1[i])**3)+G*u_n[i+1]+H*u_n[i-1]

        elif version == 'vectorized':   # (1:-1 slice style)
            u[1:-1] = A*(u_n[1:-1])+B*(u_nm1[1:-1])+C*((u_n[1:-1])**3)+D*((u_n[1:-1])**2)*(u_nm1[1:-1])+E*(u_n[1:-1])*((u_nm1[1:-1])**2)+\
            F*((u_nm1[1:-1])**3)+G*(u_n[2:])+H*(u_n[0:-2])

        elif version == 'vectorized2':  # (1:Nx slice style)
            u[1:Nx] = A*(u_n[1:Nx])+B*(u_nm1[1:Nx])+C*((u_n[1:Nx])**3)+D*((u_n[1:Nx])**2)*(u_nm1[1:Nx])+E*(u_n[1:Nx])*((u_nm1[1:Nx])**2)+\
            F*((u_nm1[1:-1])**3)+G*(u_n[2:Nx+1])+H*(u_n[0:Nx-1])

        # Insert boundary conditions
        i = Ix[0]
        if U_0 is None:
            # Set boundary values
            # x=0: i-1 -> i+1 since u[i-1]=u[i+1] when du/dn=0
            # x=L: i+1 -> i-1 since u[i+1]=u[i-1] when du/dn=0
            ip1 = i+1
            im1 = ip1
            u[i] = A*u_n[i]+B*u_nm1[i]+C*((u_n[i])**3)+D*((u_n[i])**2)*(u_nm1[i])+E*(u_n[i])*((u_nm1[i])**2)+\
            F*((u_nm1[i])**3)+G*u_n[ip1]+H*u_n[im1]
        else:
            u[0] = U_0(t[n+1])

        i = Ix[-1]
        if U_L is None:
            im1 = i-1
            ip1 = im1
            u[i] = A*u_n[i]+B*u_nm1[i]+C*((u_n[i])**3)+D*((u_n[i])**2)*(u_nm1[i])+E*(u_n[i])*((u_nm1[i])**2)+\
            F*((u_nm1[i])**3)+G*u_n[ip1]+H*u_n[im1]
        else:
            u[i] = U_L(t[n+1])

        if user_action is not None:
            if user_action(u, x, t, n+1):
                break

        # Update data structures for next step
        #u_nm1[:] = u_n;  u_n[:] = u  # safe, but slower
        u_nm1, u_n, u = u_n, u, u_nm1

    # Important to correct the mathematically wrong u=u_nm1 above
    # before returning u
    u = u_n
    cpu_time = time.perf_counter() - t0
    return u, x, t, cpu_time

###

def viz(
    omega, dt, dx, mu, L, T,  # PDE parameters
    umin, umax,               # Interval for u in plots
    U_0, U_L,
    animate=True,             # Simulation with animation?
    solver_function=solver   # Function with numerical algorithm
):
    """Run solver and visualize u at each time level."""
    import time, glob, os

    class PlotMatplotlib:
        def __call__(self, u, x, t, n):
            """user_action function for solver."""
            if n == 0:
                plt.ion()
                self.lines = plt.plot(x, u, 'r-')
                plt.xlabel('x')
                plt.ylabel('u')
                plt.axis([0, L, umin, umax])
                plt.legend(['t=%f' % t[n]], loc='lower left')
            else:
                self.lines[0].set_ydata(u)
                plt.legend(['t=%f' % t[n]], loc='lower left')
                plt.draw()
            time.sleep(2) if t[n] == 0 else time.sleep(0.2)
            plt.savefig('/Users/carlestapi/Desktop/vdp_%04d.png' % n)  # movie

    plot_u = PlotMatplotlib()

   # Clean up old movie frames
    for filename in glob.glob('tmp_*.png'):
       os.remove(filename)

    # Call solver and do the simulaton
    user_action = plot_u if animate else None
    u, x, t, cpu = solver_function(
        omega, dt, dx, mu, L, T, U_0 , U_L, user_action)

    return cpu

###

if __name__ == '__main__':
viz(
    omega, dt, dx, mu, L, T,  # PDE parameters
    umin, umax,               # Interval for u in plots
    animate=True,             # Simulation with animation?
    tool='matplotlib',        # 'matplotlib' or 'scitools'
    solver_function=solver    # Function with numerical algorithm
)


\end{lstlisting}

\bibliographystyle{unsrt}
\bibliography{references}

\begin{thebibliography}{10}

\bibitem{pikovsky2007synchronization}
Arkady Pikovsky and Michael Rosenblum.
\newblock Synchronization.
\newblock {\em Scholarpedia}, 2(12):1459, 2007.

\bibitem{bailey1988kaempfer}
Beatrice M~Bodart Bailey.
\newblock Kaempfer restor'd.
\newblock {\em Monumenta Nipponica}, pages 1--33, 1988.

\bibitem{buck1968mechanism}
John Buck and Elisabeth Buck.
\newblock Mechanism of rhythmic synchronous flashing of fireflies.
\newblock {\em Science}, 159(3821):1319--1327, 1968.

\bibitem{strogatz1993coupled}
Steven~H Strogatz, Ian Stewart, et~al.
\newblock Coupled oscillators and biological synchronization.
\newblock {\em Scientific American}, 269(6):102--109, 1993.

\bibitem{danino2010synchronized}
Tal Danino, Octavio Mondrag{\'o}n-Palomino, Lev Tsimring, and Jeff Hasty.
\newblock A synchronized quorum of genetic clocks.
\newblock {\em Nature}, 463(7279):326--330, 2010.

\bibitem{prindle2012sensing}
Arthur Prindle, Phillip Samayoa, Ivan Razinkov, Tal Danino, Lev~S Tsimring, and
  Jeff Hasty.
\newblock A sensing array of radically coupled genetic/biopixels/'.
\newblock {\em Nature}, 481(7379):39--44, 2012.

\bibitem{mcmillen2002synchronizing}
David McMillen, Nancy Kopell, Jeff Hasty, and JJ~Collins.
\newblock Synchronizing genetic relaxation oscillators by intercell signaling.
\newblock {\em Proceedings of the National Academy of Sciences},
  99(2):679--684, 2002.

\bibitem{winfree1967biological}
Arthur~T Winfree.
\newblock Biological rhythms and the behavior of populations of coupled
  oscillators.
\newblock {\em Journal of theoretical biology}, 16(1):15--42, 1967.

\bibitem{1985kuramoto}
Kuramoto.
\newblock Y. kuramoto: Chemical oscillations, waves, and turbulence,
  springer-verlag, berlin and new york, 1984, viii+ 156 , 25$\times$ 17cm,
  9,480 (springer series in synergetics, vol. 19).
\newblock {\em Springer Series in Synergetics}, 40(10):817--818, 1985.

\bibitem{hong2002synchronization}
Hyunsuk Hong, Moo-Young Choi, and Beom~Jun Kim.
\newblock Synchronization on small-world networks.
\newblock {\em Physical Review E}, 65(2):026139, 2002.

\bibitem{hong2005collective}
H~Hong, Hyunggyu Park, and MY~Choi.
\newblock Collective synchronization in spatially extended systems of coupled
  oscillators with random frequencies.
\newblock {\em Physical Review E}, 72(3):036217, 2005.

\bibitem{pikovsky2015dynamics}
Arkady Pikovsky and Michael Rosenblum.
\newblock Dynamics of globally coupled oscillators: Progress and perspectives.
\newblock {\em Chaos: An Interdisciplinary Journal of Nonlinear Science},
  25(9):097616, 2015.

\bibitem{hong2004collective}
Hyunsuk Hong, Hyunggyu Park, and Moo~Young Choi.
\newblock Collective phase synchronization in locally coupled limit-cycle
  oscillators.
\newblock {\em Physical Review E}, 70(4):045204, 2004.

\bibitem{daido1988lower}
Hiroaki Daido.
\newblock Lower critical dimension for populations of oscillators with randomly
  distributed frequencies: a renormalization-group analysis.
\newblock {\em Physical review letters}, 61(2):231, 1988.

\bibitem{bahiana1994order}
M~Bahiana and MSO Massunaga.
\newblock Order in two-dimensional oscillator lattices.
\newblock {\em Physical Review E}, 49(5):R3558, 1994.

\bibitem{aoyagi1991frequency}
Toshio Aoyagi and Yoshiki Kuramoto.
\newblock Frequency order and wave patterns of mutual entrainment in
  two-dimensional oscillator lattices.
\newblock {\em Physics Letters A}, 155(6):410--414, 1991.

\bibitem{mirollo1990amplitude}
Renato~E Mirollo and Steven~H Strogatz.
\newblock Amplitude death in an array of limit-cycle oscillators.
\newblock {\em Journal of Statistical Physics}, 60(1):245--262, 1990.

\bibitem{li2012quorum}
Zhi Li and Satish~K Nair.
\newblock Quorum sensing: how bacteria can coordinate activity and synchronize
  their response to external signals?
\newblock {\em Protein Science}, 21(10):1403--1417, 2012.

\bibitem{lee2013quantum}
Tony~E Lee and HR~Sadeghpour.
\newblock Quantum synchronization of quantum van der pol oscillators with
  trapped ions.
\newblock {\em Physical review letters}, 111(23):234101, 2013.

\bibitem{haken1985theoretical}
Hermann Haken, JA~Scott Kelso, and Heinz Bunz.
\newblock A theoretical model of phase transitions in human hand movements.
\newblock {\em Biological cybernetics}, 51(5):347--356, 1985.

\bibitem{franci2018synchronization}
Alessio Franci, Marco~Arieli Herrera-Valdez, Miguel Lara-Aparicio, and Pablo
  Padilla-Longoria.
\newblock Synchronization, oscillator death, and frequency modulation in a
  class of biologically inspired coupled oscillators.
\newblock {\em Frontiers in Applied Mathematics and Statistics}, 4:51, 2018.

\bibitem{langtangen1999computational}
Hans~Petter Langtangen.
\newblock {\em Computational partial differential equations: numerical methods
  and diffpack programming}, volume~2.
\newblock Springer Berlin, 1999.

\bibitem{garcia2004modeling}
Jordi Garcia-Ojalvo, Michael~B Elowitz, and Steven~H Strogatz.
\newblock Modeling a synthetic multicellular clock: repressilators coupled by
  quorum sensing.
\newblock {\em Proceedings of the National Academy of Sciences},
  101(30):10955--10960, 2004.

\bibitem{stricker2008fast}
Jesse Stricker, Scott Cookson, Matthew~R Bennett, William~H Mather, Lev~S
  Tsimring, and Jeff Hasty.
\newblock A fast, robust and tunable synthetic gene oscillator.
\newblock {\em Nature}, 456(7221):516--519, 2008.

\bibitem{kuwayama2013biological}
Hidekazu Kuwayama and Shuji Ishida.
\newblock Biological soliton in multicellular movement.
\newblock {\em Scientific Reports}, 3(1):1--5, 2013.

\end{thebibliography}

\end{document}